\documentclass[review]{elsarticle}

\usepackage{amsfonts,amssymb,amsmath,amsbsy,amsthm,amscd}
\usepackage{lineno,hyperref}
\usepackage[english]{babel}
\usepackage{graphicx}
\usepackage{color}
\usepackage[cp1251]{inputenc}
\usepackage{graphicx}			
\usepackage{wrapfig}			

\usepackage{geometry}			
\usepackage{indentfirst}		
\usepackage{multicol}

\usepackage{graphicx}
\usepackage{color}
\usepackage{array}

\UseRawInputEncoding

\def\dfrac#1#2{\displaystyle{#1\over #2}}

\def\bV{{\bf V}}

\def\Div{\mbox{div}\,}
\def\Rot{\mbox{rot}\,}

\def\bB{{\bf B}}

\def\bE{{\bf E}}

\newtheorem{theorem}{Theorem}[section]

\newtheorem{proposition}{Proposition}

\renewcommand{\le}{\leqslant}
\renewcommand{\ge}{\geqslant}

\begin{document}

\begin{frontmatter}

\title
{The simplest solutions of cold plasma equations: change in properties from a hydrodynamic to a  kinetic model
}

\author[1]{Olga S. Rozanova} 
\ead{rozanova@mech.math.msu.su} 
\author[2]{Lidia V. Gargyants}
\ead{gargyants@bmstu.ru} %
\address[1]{
 Lomonosov Moscow
State University, Moscow 119991 Russia}
\address[2]{Bauman Moscow State Technical University,   Moscow 105005 Russia}
\begin{abstract}
We consider the transition from the kinetic model of Landau cold plasma to the hydrodynamic one by constructing a ``multi-speed'' moment chain in the case of one spatial variable. Closing this chain at the first step leads to the standard hydrodynamic system of cold plasma. The change in the properties of the solution when closing the chain at the second step is discussed using the example of two classes of solutions - affine in space and traveling waves, and it is shown that their properties change significantly compared to the hydrodynamic model.
 \end{abstract}


\end{frontmatter}


\numberwithin{theorem}{section}
\numberwithin{remark}{section}

\section{Introduction}

In this paper, a scheme is proposed for reducing the one-dimensional kinetic equation for the electron density $F$ in phase space, corresponding to the cold plasma model, to an infinite system of quasi-hydrodynamic equations. Part of the solution components corresponds to physical macroscopic characteristics of the medium, such as the electron density, their velocity, and electric field strength. The remaining components are auxiliary in nature. When this system is closed at the first step, when the solution contains only physical components, the simplest hydrodynamic model of cold plasma is obtained.

There is a huge literature devoted to the transition from the kinetic to the hydrodynamic model, including for the multidimensional case. We  mention, without claiming to be complete, the monographs \cite{SR}, where the Vlasov-Boltzmann equation was considered as the kinetic equation, and \cite{ASR}, where the Vlasov-Maxwell-Boltzmann equation was used as the starting one, and the references contained therein. Note that several options for compiling moment chains, i.e. systems of equations linking physical quantities, including density, velocity, pressure, etc., were previously proposed. Such chains turn out to be infinite and the problem of their closure arises. The most famous is the 13-moment Grad model \cite{Grad}, which is derived using the expansion of the solution of the kinetic equation in orthogonal Hermite polynomials with the weight in the form of a Maxwell distribution (distribution of electrons by velocities at a given temperature). We also mention the model of Oraevskii et al. \cite{Oraevskii}, consisting of 16 moments. A serious problem is that the resulting system of quasilinear equations, when closed at the first step, coincides with the hydrodynamic model of cold plasma, described by a hyperbolic system, but when closed at the third step, it loses hyperbolicity. There is a large literature attempting to address this issue, we refer to \cite{Torrilhon}, \cite{Koellermeier},  \cite{Struchtrup},   \cite{Cai},  \cite{CFL}, \cite{Bobylev}.
 Note that the problem of constructing and closing moment chains arises in various contexts \cite{Kuehn}.

We propose a method in which the system preserves hyperbolicity at any step. It is hyperbolic, but not strictly hyperbolic, as is the hydrodynamic model of cold plasma obtained by closing it at the first step. The principle of closure is based on the use of the Holder inequality. Its characteristic feature is the possibility of rewriting it in terms of high-order ``velocities'', of which the first-order velocity is the usual physical velocity, and the others play an auxiliary role. The system can also be written in a conservative form, which makes it possible to study the solution numerically.

The kinetic equation in the form in which we consider it, apparently first appeared in  \cite{Landau}. Iordanskii \cite{Iord} considered the Cauchy problem for such an equation and showed that under some fairly general assumptions about the electron density
 $F(t,x,v)$ the solution preserves global smoothness for all initial data. For the hydrodynamic model of cold plasma, on the contrary, there is a wide class of initial data for which the derivatives of solution go to infinity  in a finite time.  The conditions on the initial data for the velocity and electric field strength that distinguish such a class can be found exactly \cite{RCh}.
By studying the properties of solutions of the system of quasi-hydrodynamic equations, we hope to advance in explaining this phenomenon.
It would be natural to expect that with closure at each subsequent step the requirements on the initial data that ensure the global smoothness of the solution are weakened.

In this paper, only the initial stage of this work is presented. Namely, the closure at the second step is investigated and two classes of solutions are considered: solutions linear in the spatial variable (affine solutions) and  the traveling wave. When closing the system at the second step, in addition to the main components of the solution corresponding to the electric field strength and physical velocity, we must consider an additional ``second velocity''.
Our main task is to study how the appearance of additional components of the solution affects the behavior of the main components.
Using specific classes of solutions as an example, we show that this influence is very specific. Namely, for affine solutions, the appearance of the ``second velocity'' does indeed prevent the derivatives of the main components of the solution from becoming infinite (the influence of the ``second velocity'' is similar to linear dumping), but at certain points these derivatives become discontinuous. Moreover, the solution can be continued to any point in time, but not uniquely.

Traveling wave solutions exist for both first and second step closures, but their properties change completely. For example, if for the closure at the first step a traveling wave smooth in both components exists only at sufficiently high speeds of its movement, for the second closure at any speed of movement the physical speed is smooth, and the electric field strength contains singularities. The shape of the traveling wave also changes.

 Note that the results of this work imply that the addition of the ``second speed'' should be considered as a singular perturbation. In other words, we cannot approximate the original system by directing the ``second velocity'' to zero. Closure at the second step changes the qualitative properties of the solution in comparison with the first step and leads to a system with new unusual properties. Apparently, the behavior of the solution for arbitrary initial data can only be obtained numerically. This is included in our future plans.

The paper is organized as follows. Section \ref{S2} presents kinetic and hydrodynamic models of one-dimensional cold plasma and compares the results on the smoothness of the solution to the Cauchy problem in each case. In Section \ref{S3}, equations for moments are derived and their integral characteristics, such as total mass and energy, are studied. Section \ref{S4} considers the closure of the equations at the first and second steps, shows that the first closure coincides with the hydrodynamic model of the cold plasma, and presents various forms of the equations for the second closure. Section~\ref{S5} analyzes affine solutions for the first and second closures. In Section \ref{S6} we construct traveling waves for the first and second closures. Section \ref{S7} summarizes the results obtained.

\section{Two models of cold plasma} \label{S2}

\subsection{Kinetic equation of cold plasma}

We consider electron plasma with a background of positive immobile
ions without a magnetic field in the one-dimensional case. The kinetic
equation of electron plasma without collisions for the electron density $F(t,x,v)>0$ in phase space in the one-dimensional case has the form \cite{Landau}
\begin{equation}\label{Io1}
F_t +v F_x - \frac{|e|}{m} E(t,x) F_v=0,
\end{equation}
where the electric field strength $E$, in turn, is determined from the Poisson equation
\begin{equation}\label{Io2}
 E_x= -4 \pi |e| \left(\int\limits_{\mathbb R} F dv-n_0\right),
\end{equation}
where $e$, $m$ are constant charge and mass of an electron, $n_0$ is constant ion density. In what follows we  write the equation in dimensionless form, formally setting $|e|=m=n_0=1$.

Let $n(t,x)=\int\limits_{\mathbb R} F dv$ be the electron density in physical space.
Rewrite \eqref{Io1}, \eqref{Io2} as
 \begin{equation}\label{Io}
F_t +v F_x - E(t,x) F_v=0, \qquad E_x= 1-n.
\end{equation}

Consider the Cauchy data
  \begin{equation}\label{IoCD}
F|_{t=0}= F_0(x,v)>0
\end{equation}
and the boundary condition $E\to 0$ as $x\to -\infty$.

Iordansky \cite{Iord} proved the following theorem.
\begin{theorem}
If a continuous function $F_0 (x, v)>0$ is such that
$$
\int\limits_{\mathbb R} \left(\int\limits_{\mathbb R} F_0(x,v) dv - 1 \right) dx=0,
$$
$$
F_0(x,v)<K(|v|), \qquad \int\limits_{\mathbb R} v^2 K (|v|) dv <\infty,
$$
where $K (\xi)$ is a monotonically decreasing function, then the solution of the Cauchy problem
\eqref{Io}, \eqref{IoCD} exists and is unique in the class of continuous functions $E(t,x)$ that have bounded derivative $E_x(t,x)$ on the entire
$x$-axis and vanish as $x \to \pm\infty$.
\end{theorem}

\subsection{Hydrodynamic model of cold plasma}

In vector form, the system of hydrodynamics of the electron fluid together with Maxwell's equations are
\begin{equation}
\label{base1}
\begin{array}{l}
 n_t + \Div(n \bV)=0\,,\quad
\bV_ t + \left( \bV \cdot \nabla \right) \bV
=\dfrac {e}{m} \, \left( \bE + \dfrac{1}{c} \left[\bV \times  \bB\right]\right),\vspace{0.5em}\\
\dfrac1{c}  \bE_t = - \dfrac{4 \pi}{c} e n \bV
 + {\rm rot}\, \bB\,,\quad
\dfrac1{c}  \bB_ t  =
 - {\rm rot}\, \bE\,, \quad \Div \bB=0\,,
\end{array}
\end{equation}
where $ c $ is the speed of light,
$ n, \bV $ are the density and velocity of electrons,
$ \bE, \bB $ are the vectors of the electric and magnetic fields, $x\in{\mathbb R}^3,$ $t\ge 0$, $\nabla$, $\rm div$, $\Rot $ are the gradient, divergence and vorticity with respect to the spatial variables.

In the one-dimensional case, $\bV=(V,0,0)$, $\bE=(E,0,0)$, $\bB\equiv 0$. In dimensionless form, the system \eqref{base1} can be rewritten as
 \begin{equation}
\begin{array}{c}
n_t +
\left(n\, V \right)_x
=0,\quad
V_t +
V  V_x =  - E, \quad
E_t = n\, V.
\end{array}
\label{3gl3}
\end{equation}
Suppose the solution is smooth. Then from the first and last equations (\ref{3gl3}) imply
$
\left(n +
 E_x \right)_t = 0.
$
For background density $n_0\equiv 1$ with boundary condition $E\to 0$ at $x\to -\infty$ we obtain
\begin{equation}
 n = 1 -E_x.
\label{Kn}
\end{equation}
This allows us to obtain a hyperbolic system for two components of the velocity $V$ and the electric field $E$ in the form
\begin{equation}\label{KK1}
V_t+VV_x=-E, \quad E_t+VE_x=V,
\end{equation}
where $(V,E)=(V(t,x), E(t,x))$, $t\in {\mathbb R}_+ $, $x\in {\mathbb R} $. The density $n(t,x)>0$ is found from \eqref{Kn}.
The details can be found in \cite{RCh}.

\bigskip

For \eqref{KK1} we consider the Cauchy data
\begin{equation}\label{K2}
(V,E)|_{t=0}=(V_0(x), E_0(x)).
\end{equation}
If the initial data are $C^1$ - smooth functions, then locally in $t$ there exists a smooth solution \eqref{KK1}, \eqref{K2}. However, it is known that the derivatives of the solution of such a Cauchy problem can go to infinity in a finite time, which corresponds to the formation of a shock wave, the criterion for the formation of a singularity is known \cite{RCh}: the solution preserves smoothness if and only if for any point $x_0\in\mathbb R$
\begin{eqnarray} \label {crit2}
\left (V'_0 (x_0) \right) ^ 2 + 2 \, E'_0 (x_0) -1 <0.
\end{eqnarray}

Thus, for the Cauchy problem in the kinetic formulation \eqref{Io2}, \eqref{IoCD} a continuous solution exists for any smooth initial data
$E'_0(x_0)<1$
(a natural constraint associated with $n_0>0$), whereas
for the problem in the hydrodynamic formulation \eqref{KK1}, \eqref{K2} there is no continuous solution for smooth initial data if  $E'_0(x_0)>\frac12$ at least at one point $x_0\in\mathbb R$.

Thus, a problem arises: how to move from the kinetic formulation to the hydrodynamic one and how to explain why the class of initial data that ensures the smoothness of the solution to the Cauchy problem narrows during the transition.

\section{Construction of moment equations based on the kinetic model}\label{S3}

We introduce the moments of the distribution, as is customary in probability theory, that is,
\begin{equation*}\label{Mom}
M_k(t,x)= \int\limits_{\mathbb R} v^k F(t,x,v) dv, \qquad k\in 0\cup \mathbb N,
\end{equation*}
and assume that moments of any order exist (then we can construct a chain of any length). Obviously, $n=M_0$.
Note that from the moments we can {\it formally} reconstruct $F$ by the formula
\begin{equation*}\label{FM}
F(t,x,v)=\sum\limits_{k=0}^\infty \, \frac{(-1)^k}{k!} \,\delta^{(k)}(v) \, M_k(t,x),
\end{equation*}
if $F(t,x,v)$ is understood as a generalized function with respect to the variable $v$ and its action on the test function is understood as
\begin{eqnarray*}
<F(t,x,v), \phi(v)>= \sum\limits_{k=0}^\infty \frac{(-1)^k}{k!} \, M_k(t,x) \,<\delta^{(k)}(v), \phi(v)>= \sum\limits_{k=0}^\infty \frac{1}{k!} \, M_k(t,x) \, \phi^{(k)}(v), \quad \phi(v)\in C_0^\infty({\mathbb R}).
\end{eqnarray*}
We introduce the velocity of order $k$ using the following formula:
\begin{equation}\label{Uk}
U_k(t,x) =\frac{\int\limits_{\mathbb R} v^k F(t,x,v) dv}{\int\limits_{\mathbb R} v^{k-1} F(t,x,v) dv}, \qquad k\in \mathbb N.
\end{equation}
If
\begin{equation}\label{F}
F(t,x,v)=\delta(v-\mathcal{V}(t,x)) P(t,x),
\end{equation}
where $\mathcal{V}(t,x)$ is some function, for example, $U_1$, then
$U_k=U_{k+1},k\in \mathbb N, $
otherwise it is not so.

From the first equation \eqref{Io} it follows (after multiplying by $v^k$ and integrating over $v$):
\begin{equation}\label{1}
(M_k)_t + (M_{k+1})_x = - k E M_{k-1}, k\in \mathbb N,
\end{equation}
for $k=0$ we obtain the continuity equation
\begin{equation}\label{2}
(M_0)_t + (M_{1})_x =0, \qquad (M_1 = M_0 U_1).
\end{equation}

From this we see that the physical velocity (involved in the continuity equation) is $U_1$.
From the Holder inequality (for odd $k$) it follows
\begin{equation}\label{Hold}
(M_k)^2\le M_{k+1} M_{k-1},
\end{equation}
under assumption \eqref{F} this inequality becomes an equality.
We denote
\begin{equation*}\label{lam}
\lambda_k(t,x)= ( M^2_k/M_{k-1}-M_{k+1})_x,
\end{equation*}
then  equations \eqref{1} can be rewritten as
\begin{equation*}\label{3}
(M_k)_t + ( M^2_k/M_{k-1})_x = -k E M_{k-1}+ \lambda_k, k\in \mathbb N.
\end{equation*}

{\it The principle of cutting off the chain and closing the system is that at step $k$ we assume $\lambda_k=0$.} Since $F$ decreases rapidly as $|v|\to\infty$ due to the assumption that all moments exist, then $\lambda_k$ decrease with $k$, and this principle makes sense.
The closure of the system \eqref{1}, \eqref{2} at step $k$ has the form
\begin{equation}\label{4}
M_t + A(M) M_x= B, \quad M=(M_0,..., M_k)^T, \quad B=(0,-EM_0, -2 EM_1,...,-k E M_{k-1})^T.
\end{equation}
Matrix $A$ consists of $k-1$ rows having 1 at $i+1$ place in $i$-th row and zeros at the remaining places, $i=1,...,k-1$,
and the last, $k$ -th row, consisting of zeros, except for the elements $A_{k \, k-1}= -\frac{M_k^2}{M_{k-1^2}}$ and $A_{k \, k}= \frac{2 M_k}{M_{k-1}}$.
For example, when closing at step 2
\begin{eqnarray*}
A=\begin{pmatrix}
0 & 1 & 0\\
0 & 0 & 1\\
0 & -\frac{M_2^2}{M_{1}^2} &\frac{2 M_2}{M_{1}}
\end{pmatrix}.
\end{eqnarray*}

The eigenvalues of the matrix $A=A(M_1,..., M_k)$ are real. Among them, there are two coincident values $M_k/M_{k-1}=U_k$, and the rest are zeros.
This is a non-strictly hyperbolic system, there is no complete set of eigenvectors.

System \eqref{4} can be rewritten in conservative form as
\begin{eqnarray}\nonumber
&&(M_0)_t+(U_1 M_0)_x=0,\\
&&(M_0 U_1)_t+(U_1 U_2 M_0)_x=- E M_0, \nonumber\\
&&(M_0 U_1 U_2)_t+(U_1 U_2 U_3 M_0)_x=- 2 E M_0 U_1,\label{5}\\
&&\dots\nonumber\\
&&(M_0 U_1 U_2\cdots U_k)_t+(U_1 U_2 U_3\cdots U_k^2 M_0)_x=- k E M_0 U_1 \cdots U_{k-1}.\nonumber
\end{eqnarray}
Note that if $U_1=U_2=...=U_k$, then the third and subsequent equations are consequences of the first two.

System \eqref{5} can be supplemented by the equation
\begin{equation}\label{EEE}
E_t + U_1 E_x= U_1,
\end{equation}
which is a consequence of the first equation \eqref{5} and the second equation \eqref{Io}.

  If we introduce a new variable ${\mathcal E}=E M_0$, then \eqref{EEE} can be rewtitten in a conservative form as
  \begin{equation*}\label{EEEE}
 {\mathcal E}_t + (U_1 {\mathcal E})_x= M_0 U_1.
\end{equation*}

A significant problem for the obtained model is the possibility of $M_k$ vanishing for odd $k$, since this quantity appears in the denominator of both $\lambda_k$ and the coefficients of the matrix $A$. This problem can be avoided, for example, by considering the distribution $F(t,x,v)$ only on the semiaxis $v>0$ with zero boundary condition for $v=0$. In this case, the Holder inequality  implies \eqref{Hold}, as well as
\begin{equation*}\label{UUU}
0\le U_k \le U_{k+1}, \qquad k \in \mathbb N.
\end{equation*}
Note that the inequality
\begin{equation}\label{U12}
|U_1| \le | U_{2}|
\end{equation}
takes place without any additional assumptions on the properties of $F$.
In addition, we note that from  definition  \eqref{Uk} it follows that  at those points where $U_1=0$ the quantity $U_2$ tends to infinity, if $F$ has no delta-singularity at $v=0$. The latter case corresponds to the cold plasma equations, the first step of closure. In this case $U_1=U_2=0$.

\subsection{Integral characteristics}

As is well known (e.g., \cite{Kuz_Zakh}), system \eqref{Io}  has two conserved integrals corresponding to mass and energy
\begin{equation*}\label{CLF}
\int\limits_{\mathbb R}\int\limits_{\mathbb R}  F dv dx  = {\mathbb C}={\rm const}, \quad \dfrac12\,\int\limits_{\mathbb R}\int\limits_{\mathbb R} (v^2 F + E^2) dv dx = {\mathbb E}={\rm const}.
\end{equation*}

All chains \eqref{5} also have these conserved quantities (under the assumption of sufficient decrease of the components of the solution at infinity), which in this case have the form
\begin{equation*}\label{CLM}
\int\limits_{\mathbb R} M_0 dx = {\mathbb C}={\rm const}, \quad 
\dfrac12\,\int\limits_{\mathbb R} (M_0 U_1 U_2 + E^2) dx = {\mathbb E}={\rm const}.
\end{equation*}
The conservation of mass follows from the first equation \eqref{5}, and the conservation of energy follows from the first three equations \eqref{5} and \eqref{E}, if we take into account condition \eqref{Kn}.

\section{Form of  system \eqref{5}, \eqref{EEE} for $k=1$ and $k=2$}\label{S4}

The closure at the first step has a conservative form
 \begin{eqnarray*}\label{6}
{\mathcal E}_t + (U_1 {\mathcal E})_x= M_0 U_1,\quad
  (M_0)_t+(U_1 M_0)_x=0,\quad
  (M_0 U_1)_t+\left(U^2_1 M_0\right)_x=-  {\mathcal E},
  \end{eqnarray*}
 which on smooth solutions is equivalent to the hydrodynamic equations of cold plasma \eqref{KK1}.

  The closure at the second step in a conservative form in variables $\mathcal E$, $M_0$, $U_1$, $U_2$ is
  \begin{eqnarray*}
   &{\mathcal E}_t + (U_1 {\mathcal E})_x= M_0 U_1,\quad
 & (M_0)_t+(U_1 M_0)_x=0,\label{8}\\
  &(M_0 U_1)_t+(U_1 U_2 M_0)_x=- \mathcal E,\quad & (M_0 U_1 U_2)_t+(U_1 U_2^2 M_0)_x=- 2 \mathcal E U_1.\nonumber
  \end{eqnarray*}
   This system can be rewritten  in terms of  ${\mathcal E}$, $M_0$, $M_1$, $M_2$ as
   \begin{eqnarray}\label{8.1}
 &{\mathcal E}_t + \left(\dfrac{M_1}{ M_0}\, {\mathcal E}\right)_x= M_1,\quad
 & (M_0)_t+(M_1)_x=0,\\
  &(M_1)_t+(M_2)_x=-  {\mathcal E},\quad &(M_2)_t+\left(\frac{M_2^2}{ M_1}\right)_x=- 2 {\mathcal E} \,\frac{M_1}{ M_0},\nonumber
  \end{eqnarray}
and in terms $E, U_1, U_2$
 as
\begin{eqnarray}\label{9}
&& E_t + U_1 E_x= U_1, \nonumber\\
&&  (U_1)_t+(U_2-U_1) (U_1)_x +U_1 (U_2)_x=-E + (U_2-U_1) U_1 \frac{E_{xx}}{1-E_x},\\
 && (U_2)_t +U_2 (U_2)_x=- 2 E +\frac{U_2}{U_1} E.\nonumber
  \end{eqnarray}

The specificity of such systems can be demonstrated already in the simplest example of solutions independent of the spatial coordinate.
Indeed, for such solutions one of the components is separated from the previous ones, which significantly simplifies the analysis. For example, for the system \eqref{9} we have $E=C_1 \cos(t+\theta)$, $U_1=-C_1 \sin(t+\theta)$, $C_1, \theta$ are constants, and $U_2= \frac{1}{\sin(t+\theta)}\left(C_2+ \frac{C_1}{2} \cos (2(t+\theta))\right)$, $C_2$ is a constant. Thus, $U_2$  becomes unbounded whenever $U_1$ vanishes, except for $U_1=U_2$.

On the other hand, it follows from \eqref{8.1} that $M_0=M_0(0)=\rm const$, ${\mathcal E}=C_1 \cos(t+\theta)$, $M_1=-C_1 \sin(t+\theta)$,
$M_2=\frac{C_1^2}{2 M_0} \cos 2(t+\theta)+C_2$, and the components of the solution no longer contain singularities.

\section{Affine solutions. Closure at step two.}\label{S5}

Recall that a solution is called affine if it is linear in the spatial coordinate $x$. We
will consider the system~\eqref{9} and seek a solution of the form
\begin{equation*}
E(t,x)=a(t)x+b(t),\quad  U_1(t,x)=\gamma_1(t)x+\delta_{10}(t),\quad U_2(t,x)=\gamma_2(t)x+\delta_{20}(t).
\end{equation*}
It is easy to verify that in the general case, substituting such expressions into \eqref{9} leads to an overdetermined system of ordinary differential equations. In order to remove the overdeterminedness, we have to assume that
$\dfrac{\gamma_1}{\gamma_2}=\dfrac{\delta_{10}}{\delta_{20}}=\rm const$, which means that the expression $U_2/U_1$ does not contain a singularity.
Next, we can shift the origin to the point $\bar x=\frac{\delta_{10}}{\gamma_1}$ and consider the components of the solution $U_1$ and $U_2$ as not containing zeroth-order terms in $x$ from the very beginning. Calculations show that the remaining component $E$ also does not contain zeroth-order terms in $x$. Thus, to determine the three unknown functions, we obtain the system
\begin{eqnarray}\label{agg}
\dot a=-\gamma_1(a-1),\quad
\dot\gamma_1=-a-2\gamma_1\gamma_2+\gamma_1^2,\quad
\dot\gamma_2=-2a+a\dfrac{\gamma_2}{\gamma_1}-\gamma_2^2,
\end{eqnarray}
which we consider together with the initial conditions
\begin{eqnarray*}\label{aggCD}
a(0)=a_0, \quad \gamma_1(0)=\gamma_{10}, \quad \gamma_2(0)=\gamma_{20}.
\end{eqnarray*}
Note that from the condition $n>0$ it follows that $a<1$.

In addition, it follows from \eqref{U12} that $|\gamma_1|\le |\gamma_2|$. For $\gamma_1= \gamma_2$ we obtain an affine solution, for which $U_1=U_2$, which corresponds to the first closure of the chain, that is, the original hydrodynamic model of the cold plasma \eqref{KK1}. From \eqref{agg} we have
\begin{eqnarray}\label{ag}
\dot a=-\gamma_1(a-1),\quad
\dot\gamma_1=-a-\gamma_1^2.
\end{eqnarray}
This system has a first integral
\begin{eqnarray*}\label{intag}
\gamma_1^2+2a-1 +C(a-1)^2=0,
\end{eqnarray*}
from which we immediately obtain that the phase trajectories of the equation \eqref{ag} are bounded and the derivatives of the solution $(E, V)$ (i.e. $a$ and $\gamma_1$) are bounded for all $t$
if and only if
\begin{eqnarray}\label{critag}
\gamma_{10}^2+2 a_0-1 < 0.
\end{eqnarray}
In this case, the motion on the plane $(a, \gamma_1)$ occurs along ellipses containing the origin of coordinates inside.
It is easy to see that for affine solutions this criterion coincides with \eqref {crit2}.
Thus, if the initial data do not satisfy  condition \eqref{critag}, then the derivatives of the solution tend to infinity during time $t^*<2\pi$. Otherwise, the derivatives remain bounded for all $t>0$.

We will show that if $\gamma_{20}>\gamma_{10}$, then the set of initial data $(a_0,\gamma_{10} )$ corresponding to the globally time smooth solution $(E,U_1)$ expands.

Namely, we will prove the following statement.

\begin{proposition}{\it
1. There exist initial data $(a_0,\gamma_{10}, \gamma_{10}<\gamma_{20} )$ such that the solution of the Cauchy problem \eqref{ag}, $(a_0,\gamma_{10})$, becomes infinite during time $t^*<2\pi$, while the components $(a,\gamma_{1})$ of the solution of  system \eqref{agg} corresponding to the same data remain bounded until time $t_1>t^*$, while the solution exists in the classical sense. For $t\to t_1-0$ the components $(a,\gamma_{1})$ have a finite limit $(a_-,0)$, and $\gamma_2$ becomes $+\infty$.

2. For these initial data, there are infinitely many solutions of  system \eqref{agg} such that the components $(a,\gamma_{1})$ are bounded for all $t>0$, but at the points where $\gamma_1$ vanishes, they have a discontinuity; the component $\gamma_2$ at these points becomes $+\infty$.}
\end{proposition}

\proof

By making the substitution $\varepsilon=\gamma_2-\gamma_1$,
instead of the last two equations of the system \eqref{agg} we get
\begin{eqnarray*}
\dot\gamma_1=-a-\gamma_1^2-2\gamma_1\varepsilon,\quad
\dot\varepsilon=a\dfrac{\varepsilon}{\gamma_1}-\varepsilon^2,
\end{eqnarray*}
that after the replacement $\gamma_1=\dfrac{a}{q}$ turns into
\begin{eqnarray}\label{qe}
\dot q=1+2\varepsilon q +q^2,\quad \dot\varepsilon=-\varepsilon^2+q\varepsilon.
\end{eqnarray}
Note that \eqref{U12} implies the condition $\varepsilon\ge 0$.

For the case $\varepsilon= 0$, for which $U_1=U_2$, which corresponds to the first closure of the chain, i.e. the original hydrodynamic model of cold plasma.
Let us pass to the case $\varepsilon> 0$ and show that
there exists a pair
$(a_0,\gamma_{10}<0)$ such that the solution of the problem \eqref{agg} with the data $ (a_0,\gamma_{10},\gamma_{20}=\gamma_{10})$ blows up and goes to minus-infinity in a finite time $t^*<2\pi$, whereas for the solution of the same problem with the data $ (a_0,\gamma_{10},\gamma_{20}>\gamma_{10})$  remains bounded until $\gamma_{1}$ turns to zero; the solution cannot be continued as a classical one beyond the point $t_1$ at which it first occurs.

First of all, we will show that for all possible values of $(a_0, \gamma_{10}, \gamma_{20}>\gamma_{10})$, there exists $\varepsilon_*> 0$ such that $\varepsilon>\varepsilon_*> 0$ for all $t$ for which the classical solution of  \eqref{agg} exists. To do this, consider the phase trajectories of the system \eqref{qe} on the plane $(q,\varepsilon )$ and note that for $\varepsilon_0=\gamma_{20}-\gamma_{10}>0$ all trajectories lie in the half-plane $\varepsilon>0$ and have a minimum at those points where $q=\varepsilon$, increase in $q$ in those regions where $\varepsilon<q,$ $q\ge 0$, as well as $\varepsilon<-\dfrac{1+q^2}{2q}$, $q<0$. In the remaining regions of the half-plane $\varepsilon>0$ the phase trajectories decrease in $q$. Thus, for any initial data there is a minimum point of the phase trajectory at which the value $\varepsilon_*> 0$ is reached. By setting $(a_0, \gamma_1\ne 0,\gamma_{20}>\gamma_{10}$, we set $\varepsilon_*> 0$. Note that the limit as $q\to -\infty$ does not exist: all trajectories turn and enter the region $\varepsilon>-\dfrac{1+q^2}{2q}$. This can be shown by examining the curves on which the direction of convexity of the trajectories changes. For $q<0$, this is the curve implicitly defined as $4\varepsilon^2+q\varepsilon+6q\varepsilon^3+q^2-3 q^3\varepsilon+1+6q^2\varepsilon^2 = 0$. This curve consists of two pieces.  One is entirely in the region $\varepsilon>-\dfrac{1+q^2}{2q}$, the other is the boundary of this region.
Fig.\ref{aff} demonstrates the behavior of phase trajectories of the system \eqref{qe} for different initial data.

\begin{figure}[!h]
\begin{multicols}{2}
\hskip1cm
\includegraphics[scale=0.7]{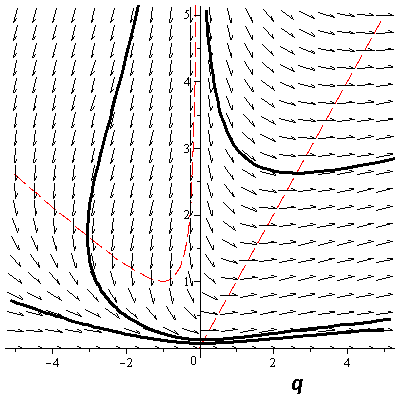}
\caption{
Behavior of phase trajectories of  system \eqref{qe}}
\label{aff}
\hskip1cm
\includegraphics[scale=0.7]{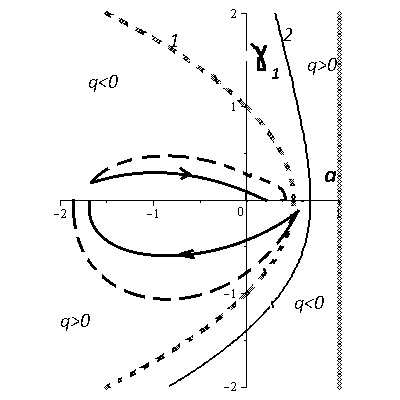}
\caption{
Continuation of the solution for $t>0$, non-uniqueness.}
\label{num1nm}
\end{multicols}
\end{figure}

Note that on the plane $(a,\gamma_1)$ for $\varepsilon=0$, as $t$ increases in the half-plane $\gamma_1<0$, the motion occurs from right to left, i.e. $q=\dfrac{a}{\gamma_1}$ changes from negative to positive values. If the projection of the phase trajectory of the system \eqref{agg} onto the plane $(a,\gamma_1)$ is bounded, then $\gamma_1$ vanishes in a finite time and in the plane $q>0$, with increasing time, the variable $q$ tends to infinity in a finite time. If $\gamma_1\to -\infty$, then $q$  tends to a constant over time.

Let us consider the first and second equations \eqref{agg}, namely
\begin{eqnarray}\label{age2}
\dot a=-\gamma_1(a-1),\quad
\dot\gamma_1=-a-\gamma_1^2-2 \gamma_1 \varepsilon.
\end{eqnarray}
Since for $a<1$ we have $\gamma_1<0$, then the following estimate is valid:
\begin{eqnarray}\label{age}
\dfrac {d \gamma_1}{da}= \dfrac{a+\gamma_1^2+2 \gamma_1 \varepsilon}{\gamma_1(a-1)}<\dfrac{a+\gamma_1^2+2 \gamma_1 \varepsilon_*}{\gamma_1(a-1)}.
\end{eqnarray}
The variable $a$ decreases for $\gamma_1<0$, then, as follows from Chaplygin's theorem on differential inequalities, the projection of the phase trajectory of  system \eqref{age2}, outgoing from the point $(a_0,0)$ is bounded from below for $\gamma_1<0$ by the trajectory of the system
\begin{eqnarray}\label{ages}
\dot a=-\gamma_1(a-1),\quad
\dot\gamma_1=-a-\gamma_1^2-2 \gamma_1 \varepsilon_*.
\end{eqnarray}
However, as shown in \cite{Rozanova:RChD_fric}, this system describes the behavior of the derivatives of the solution of the hydrodynamic equations of a cold plasma with a constant collision coefficient $\nu>0$
\begin{equation*}\label{K1nu}
V_t+VV_x=-E-\nu V, \quad E_t+VE_x=V, \quad \nu=\varepsilon_*.
\end{equation*}
For such a system, a criterion for the formation of a singularity is known. The region on the plane $(a_0,\gamma_{10})$ corresponding to globally smooth solutions is wider for $\nu>0$ than the same region for $\nu=0$. In particular, if for $\nu=0$ the point $(a_0,0)$ falls into the region of initial data corresponding to globally smooth solutions only for $a_0<\frac12$, then for $\nu>0$ the value of $a_0$ can be chosen to be greater than $\frac12$ (it increases with increasing $\nu$).
Let $\varepsilon_*>0$ be such that the phase trajectory of  system \eqref{ages} outgoing from $(a_0,\gamma_{10}<0)$ does not tend to infinity, but intersects the axis $\gamma_1=0$ at the point $(a_1^*,0)$. The existence of such a trajectory follows from the results of \cite{RCh}. Then the projection of the phase trajectory of  \eqref{age2} outgoing from the point $(a_0,0)$ intersects the axis $\gamma_1=0$ at the point $(a_1,0)$, where $a_1^*<a_1<0$. At this point, $q$ becomes $+\infty$, and $\varepsilon$ becomes positive infinity. This means that $\gamma_2$ also becomes positive infinity, and the smoothness of the component $U_2$ is lost at this point. 
Thus, this solution of  system \eqref{agg} ceases to exist as a classical one. Point 1 is proven.

\medskip

Let us prove point 2. First of all, we note that in the half-plane $\gamma_1>0$ the estimate \eqref{age} is also valid, however the variable $a$ increases and the projection of the phase trajectory of system \eqref{agg}, emerging from the point $(a_2<0,\gamma_{10}>0)$ for all possible $\gamma_{20}$ is bounded from above by the trajectory of  system \eqref{ages}, but possibly with a different value of $\varepsilon_*$.

The problem is how to find new initial conditions $(a_0<0,\gamma_{10}>0, \gamma_{20})$ when going from the lower to the upper half-plane in order to glue the new phase curve to the phase curve from the lower half-plane. The same problem arises if we then want to go from the upper to the lower half-plane again, that is, whenever $q$ becomes $+\infty$. Fig.\ref{num1nm} illustrates these considerations.
Curve 1 delimits the domain of initial data of system \eqref{ages} for $\varepsilon_* = 0$ so that for initial data from the domain containing the origin, the solution does not tend to infinity (the phase trajectory is an ellipse). Curve 2 delimits the domain of initial data in the same way for $2>\varepsilon_* >0$, the phase trajectory is a spiral.

In order to construct a bounded continuation of the solution of  \eqref{agg} with initial data $(a_0,\gamma_{10}<0,\gamma_{20}<0)$, one can follow the following algorithm.

1. We define $\varepsilon_*$ at the first step of constructing the solution, as described in point 1.

2. For this $\varepsilon_*$ we find $\alpha^*>0$ such that the point $(\alpha^*, 0)$ belongs to the boundary of the domain of initial data corresponding to globally smooth solutions of the system \eqref{ages}. Note that $\frac12<\alpha^*<1$.

3. In the upper half-plane $\gamma_1>0$
we construct the trajectory of the solution of  system \eqref{ages} corresponding to the solution of the Cauchy problem with initial data $(\beta, 0),$ $0<\beta<\alpha^*$ in reverse time, and find the point $(\alpha_1<0, 0)$ of intersection of this trajectory with the axis $\gamma_1=0$ for $a<0$.
4. In order to take the second step and continue the solution constructed in point 1 into the upper half-plane, we choose as the initial data of the system \eqref{agg} the point $(a_2^*<0, \gamma_{12}^*>0, \gamma_{22}^*> \gamma_{12}^*>0)$
from the considerations that $(a_2^*<0, \gamma_{12}^*>0)$ lie under the trajectory constructed above, and $\gamma_{22}^*$ corresponds to the value of $\varepsilon_*$ found above. For this, it is sufficient that the point $(q_2, \varepsilon_2)$, $q_2=\dfrac{a_2^*}{ \gamma_{12}^*},$ $\varepsilon_2=\gamma_{12}^*-\gamma_{22}^*$, lies on the phase trajectory of system \eqref{qe} constructed in point 1.
Thus, the projection of the trajectory of the system \eqref{agg} onto the plane $(a, \gamma_1)$ intersects the axis $\gamma_1=0$ at the point $(\alpha_2^*, 0)$, $0<\alpha_2^*<\beta$.

5. In order to continue the solution constructed in the second step to the lower half-plane $\gamma_1<0$, that is, to take the third step, we choose as the initial data
the point $(a_3<0, \gamma_{13}<0, \gamma_{23}> \gamma_{13}>0)$ so that the point $(q_3, \varepsilon_3)$, $q_3=\dfrac{a_3}{ \gamma_{13}},$ $\varepsilon_3=\gamma_{13}-\gamma_{23}$ lies on the phase trajectory of  system \eqref{qe} constructed in step 1, and
$(a_3<0, \gamma_{13}<0$ lies above the trajectory of the solution of the system \eqref{ages}, corresponding to the solution of the Cauchy problem with the initial data $(\beta, 0)$ in the direct time. In particular, as the point $(a_3<0, \gamma_{13}<0, \gamma_{23}> \gamma_{13}>0)$ we can choose the initial value $(a_0, \gamma_{10}, \gamma_{20})$, then the solution constructed at the third step will coincide with that constructed at the first, and, continuing the procedure, we will obtain a periodic solution.

Thus, the solution with bounded but discontinuous $(a,\gamma_1)$ is constructed for all $t>0$. Obviously, it can be constructed in an infinite number of ways.

The proof is complete. $\Box$

\medskip

From \eqref{age2} we see that this system can be viewed as a matrix Riccati equation with a discontinuous (unbounded) coefficient
$\varepsilon(t)$, which prevents a unique continuation of the solution after the time at which $\gamma_1$ vanishes.

\section{Traveling waves}\label{S6}

For the closure at the first step, the traveling waves, i.e. solutions of the form $E=E(\xi), U_1=U_1(\xi),$ $ \xi=x-wt, $ $w={\rm const}$
are well known \cite{AP56}, \cite{RCh}. From the equation \eqref{KK1} it follows
\begin{equation}\label{EU11}
E'=\dfrac{U_1}{U_1-w},\qquad U_1'=-\dfrac{E}{U_1-w},
\end{equation}
this system has a first integral $U_1^2+E^2=\rm const.$ Taking this into account, from \eqref{EU11} we find $U_1$ (implicitly)
\begin{eqnarray*}
\xi+C=\sqrt{\mathcal{I}_0^2-U_1^2}+w\arcsin\dfrac{U_1}{\mathcal{I}_0},\quad \mathcal{I}_0^2=E^2(0)+U_1^2(0),\quad
C=\sqrt{\mathcal{I}_0^2-U_1^2(0)}+w\arcsin\dfrac{U_1(0)}{\mathcal{I}_0}.
\end{eqnarray*}
The solution is periodic with period $T=2\pi|w|.$ Note that from \eqref{EU11} it follows that a smooth solution exists only for $w^2\ge \mathcal{I}_0^2, $ otherwise the derivatives become unbounded.

For the closure at the second step we  also look for a solution of the form $E=E(\xi), U_1=U_1(\xi), U_2=U_2(\xi), \xi=x-wt, $ $w={\rm const}$, but we will take into account that now we must limit ourselves to the region $|U_1|\le |U_2|$. For definiteness we assume that $U_1>0$ and $w>0$.

We obtain an ODE system
\begin{eqnarray*}
&&(-w+U_1)E'=U_1,\\
&&(-w U_2-U_1)U_1'+U_1U_2'=-E+(U_2-U_1)U_1\dfrac{E''}{1-E'},\\
&&(-w+U_2)U_2'=-2E+\dfrac{U_2}{U_1}E,
\end{eqnarray*}
which can be written in normal form
\begin{eqnarray}
&&E'=\dfrac{U_1}{U_1-w},\label{E}\\
&&U_1'=E\,\dfrac{(U_1-w)(2U_2-2U_1-w)}{w(U_2-w)^2},\label{U1}\\
&&U_2'=E\,\dfrac{U_2-2U_1}{U_1(U_2-w)}.\label{U2}
\end{eqnarray}

Note that from \eqref{U1},~\eqref{U2} we have the equation
\begin{equation}\label{U2(U1)}
\dfrac{dU_2}{dU_1}=\dfrac{w(U_2-2U_1)(U_2-w)}{U_1(U_1-w)(2U_2-2U_1-w)},
\end{equation}
which would give an explicit connection between $U_1$ and $U_2$ if it were integrable. Not all solutions of this equation correspond to integral curves of system \eqref{E}-\eqref{U2}. For example, the obvious solution $U_2=w$ does not correspond to any solution of \eqref{E}-\eqref{U2}.
However, another simple solution $U_2=U_1$ corresponds to the reduction to system \eqref{EU11}.

Unfortunately, the equation is not integrable and it has to be investigated only qualitatively. The equation has five singular points:
\begin{eqnarray}\label{SP}
&A_1: U_1=0, U_2=0, \quad  A_2: U_1=0, U_2=w,\\
& A_3: U_1=w, U_2=w,\quad   A_4: U_1=w, U_2=2w,\quad A_5: U_1=\frac12 w, U_2=w.\nonumber
\end{eqnarray}
Linear analysis shows that $A_1$ and $A_2$ are saddles, and $A_3$, $A_4$, $A_5$ are nodes ($A_5$ is a degenerate node).

Note that the line $U_2=U_1$ is a separatrix of the saddle $A_1$ leading to the node $A_3$. Only it is a projection of the phase curves of the original system \eqref{E}-\eqref{U2} with the constraint $|U_1|\le |U_2|$, along which motion to the half-plane $U_1<0$ through the point $A_1$ and a change in the sign of the derivative $E$ is possible. For the remaining phase trajectories, the intersection with the axis $U_1=0$  does not occur. Since the singular points \eqref{SP} are not in the lower half-plane $U_1<0$, we will further analyze the behavior of phase trajectories that satisfy the condition $U_2>U_1> 0$. Under this condition, according to \eqref{E}, $E$ increases for $U_1>w$ and decreases for
$0<U_1<w$; the projections of phase trajectories cannot intersect the line $U_1=w$ according to \eqref{U2(U1)}.

The solutions of  system \eqref{E}-\eqref{U2} are classified depending on where the projection of its phase curves is located on the plane
$(U_1,U_2)$.

\begin{figure}[!h]
\begin{multicols}{2}
\hskip-3cm
\includegraphics[scale=0.7]{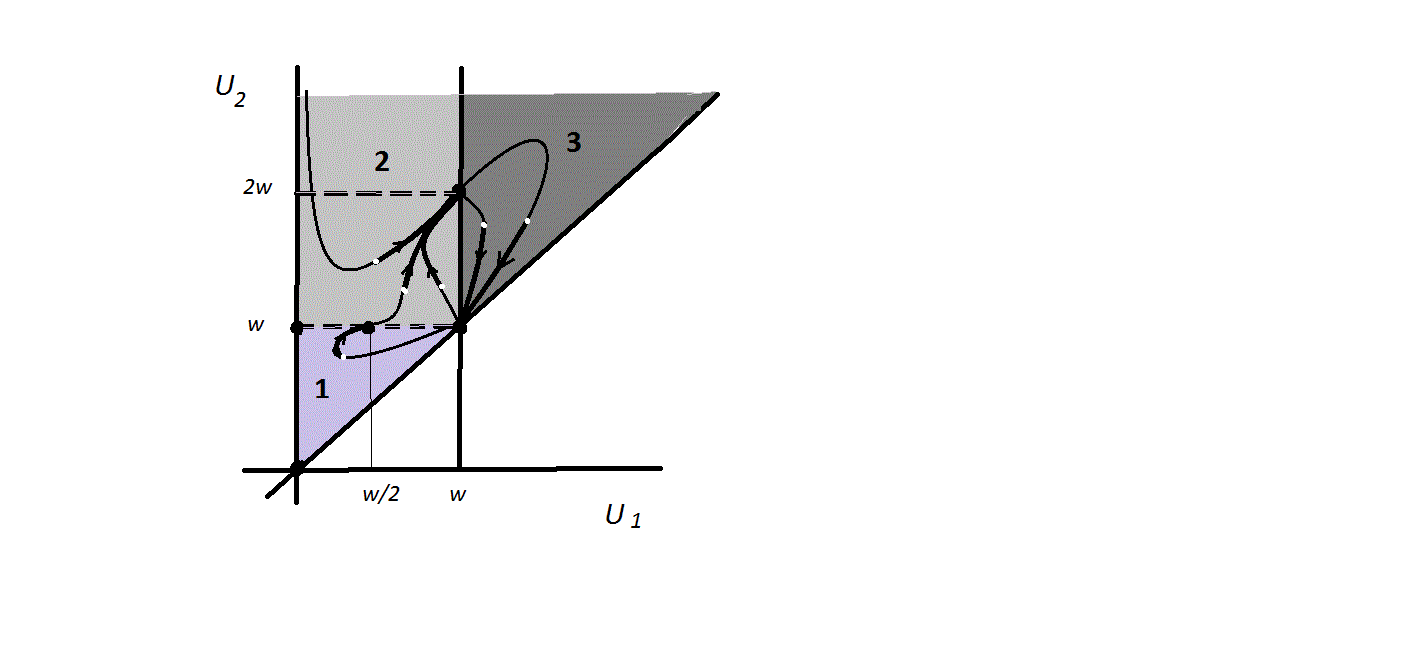}
\caption{\\
Projections of phase curves of  system \eqref{E}-\eqref{U2} on the plane $(U_1, U_2)$ depending on the region. Arrows indicate the motion with increasing $\xi$ at $E<0$, at the points marked in white the sign of $E$ changes and the direction of motion changes.}
\label{123zone}
\end{multicols}
\end{figure}

 1. Region $w>U_2>U_1> 0$. The motion starts from $A_5$ and returns to the same point in finite time, $E$ decreases, remains bounded, the direction of motion changes at $E=0$, at the end points of the trajectory the derivatives of the solution become unbounded. 

 2. Region $U_2>w>U_1> 0$. The projection of the trajectory leaves $A_4$ and returns there in a finite time, $E$ decreases, the direction of motion changes at $E=0$. At the end points of the trajectory, $E$ becomes unlimited. 

3. Region $U_2>U_1>w$. The projection of the trajectory leaves $A_3$ and arrives there in a finite time. $E$ increases, remains bounded, the direction of motion changes at $E=0$. At the end points of the trajectory $E$, the derivatives of the solution become unbounded.

Fig.\ref{123zone} shows the location of the projections of the phase curves of system \eqref{E}-\eqref{U2}  on the plane $(U_1, U_2)$.

Let us prove the statements about  the behavior of the solution constructed above.

\begin{proposition}

1. The support (the domain of definition) of a continuous solution of system \eqref{E}-\eqref{U2} is always bounded for $U_2\ne U_2$, its size depends on the choice of initial conditions.

2. At the end points of the support

a) \, $E$ is bounded in regions $1$ and $3$ and unbounded in region $2$;

b) \, $U_1$ and $U_2$ are bounded in all regions;

c) \, $E'$ is bounded in region $1$ and unbounded in regions $2$ and $3$;

d) \, $U_1'$ and $U_2'$ are bounded in region 2 and unbounded in regions 1 and 3.
\end{proposition}

\proof
Despite the fact that  equation \eqref{U2(U1)} cannot be explicitly integrated, it is possible to find the asymptotics of the solution near the singular points, that is, to approximately find the dependence $U_2(U_1)$. Then, from the system \eqref{E}-\eqref{U2}, the dependence $E(U_1)$ can be approximately found, which makes it possible to obtain from \eqref{U1} an equation connecting $U_1$ and $\xi$, and to separate the variables, obtaining an expression of the form
$ \Psi(U_1) d U_1= d \xi$. The convergence of the integral $\displaystyle\int\limits_{U^*_1}\Psi(U_1) d U_1$ at the point $U_1=U^*_1=\dfrac{w}{2}$ (for region 1) and $U_1=w$ (for regions 2 and 3)  means the boundedness of the support. In some cases it may be more convenient to use equation \eqref{U2} in a similar manner.

\medskip

 Case 1. We study the behavior of integral curves near the singular point $A_5$. Let us use the substitution $U_2(U_1)=w+\Psi(U_1-\frac{w}{2})$,
 where $\Psi(\eta)$, $\eta=U_1-\frac{w}{2}$, is a function differentiable at $\eta=0$, $\Psi(0)=0$.
 Since the right hand side of \eqref{U2(U1)}
 as $U_1\to \dfrac{w}{2}$, $U_2\to {w}$ is equivalent to $\dfrac{-4 (U_2-2U_1)(U_2-w) }{w (2U_2-2U_1-w)}$, therefore
 \begin{equation}\label{Psi}\dfrac{d \Psi}{d \eta}=-\frac{2}{w}\frac{(\Psi-2\eta) \Psi}{\Psi-\eta}=-\frac{2}{w}\frac{(1-2\frac{\eta}{\Psi}) }{1-\frac{\eta}{\Psi}}\, \Psi,\end{equation}
whence it follows that $\Psi\sim \eta $ as $\eta\to 0$. Indeed, if $\frac{\Psi}{\eta}\to 0$,   $\frac{\Psi}{\eta}\to \infty$ or  $\frac{\Psi}{\eta}\to k \ne 1$, we see that the right hand side of \eqref{Psi}
is equivalent to $k \Psi$ with a constant $k$. Thus, $\Psi \sim C e^{k \eta}$ and cannot be zero unless $C=0$. Note that we assume that the dependence $\Psi(\eta)$ can be presented as a power series, we get
$\Psi(\eta)=\eta +\frac{2}{w}\eta^2 + o(\eta^2)$. Further,
from \eqref{E}, \eqref{U1} it follows
\begin{equation}\label{EU1}
E {dE} =\frac{w U_1 (U_2-w)^2}{ (2U_2-2U_1-w)(U_1-w)^2}{dU_1},
\end{equation}
whence we obtain that near $A_5$ the right-hand side of this equation is equivalent to
\begin{equation*}
\frac{1}{w}\,\frac{\Psi^2(\eta)}{\Psi(\eta)-\eta}\,{d\eta}, \quad \eta\to 0.
\end{equation*}
Since $\Psi(\eta)-\eta=o(\eta),\quad \eta\to 0$, then $\displaystyle \int\limits_{0}^{\epsilon} \frac{\Psi^2(\eta)}{\Psi(\eta)-\eta}\,{d\eta}<\infty$. Next,
\eqref{EU1} implies $E\to E_*= \rm const$ for $U_1\to \frac{w}{2}$. Further,  near $A_5$ the right-hand side of \eqref{E} is equivalent to a constant, which means the convergence of the integral $\displaystyle\int\limits_{\xi_*}\dfrac{U_1}{U_1-w} d \xi$, $E(\xi_*)=E_*$ and the boundedness of the carrier of the traveling wave. The boundedness of $E'$ and the unboundedness of $U_1'$, $U_2'$ follows directly from \eqref{E}-\eqref{U2}.

\begin{figure}[htb]
 \hspace{1.5cm}
\begin{minipage}{0.4\columnwidth}
\includegraphics[scale=0.6]{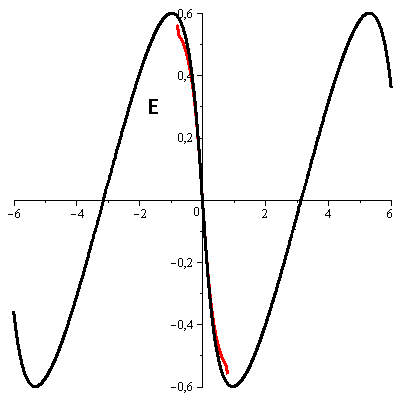}
\end{minipage}
\hspace{1.5cm}
\begin{minipage}{0.4\columnwidth}
\includegraphics[scale=0.6]{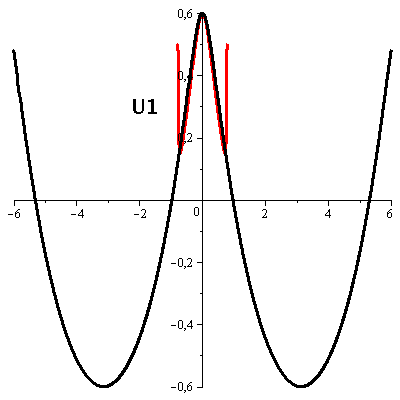}
\end{minipage}
\caption{Traveling wave for the second closure of parameter values from region 1 (red line), compared to the traveling wave for the first closure (black line), on the left  $E$, on the right $U_1$. Everywhere $w=1$, $E(0)=0$, $U_1(0)=0.6$. For the first closure $U_2(0)=0.6$, for the second closure $U_2(0)=0.7$.}\label{1EU1}
\end{figure}

Fig.\ref{1EU1} illustrates the form of the traveling wave in domain 1 in comparison with the traveling wave for $U_1=U_2$.

\medskip

Case 2. We investigate the behavior of integral curves near the singular point $A_4$.
Making the substitution $\tilde U_1=\ln(U_1-w),$ we obtain
\begin{equation*}
  \dfrac{dU_2}{d\tilde U_1}=\dfrac{w(U_2-2w-2e^{\tilde U_1})(U_2-w)}{(w+e^{\tilde U_1})(2U_2-3w-2e^{\tilde U_1})}.
\end{equation*}
For $U_1\to w$ we have $e^{\tilde U_1}\to 0,$ therefore the right-hand side of this equation near $A_4$ is equivalent to
$
\dfrac{(U_2-2w)(U_2-w)}{2U_2-3w}.
$
Thus, after integration we obtain that near $A_4$
\begin{equation}\label{SpP}
U_1\sim w+C_1 (U_2-w)(U_2-2w), \quad C_1=\rm const
\end{equation}
(the constants are not related to those mentioned in the previous paragraph, although they are named with the same letters).
Next, from equations~\eqref{E},~\eqref{U2} we obtain
\begin{equation}\label{EU2}
\dfrac{dE}{dU_2}=\dfrac{U_1^2(U_2-w)}{E(U_2-2U_1)(U_1-w)}.
\end{equation}
From ~\eqref{SpP} it follows that the right-hand side of this equation near $A_4$ is equivalent to
$$
\dfrac{(w+C(U_2-w)(U_2-2w))^2 }{C(U_2-2w)^2(1-2C(U_2-w))E}.
$$
We see that for $U_2\to 2w$ we have $E^2 \sim \dfrac{\tilde C_1}{U_2-2w}+ C_2$, $\tilde C_1(C_1), C_2$ are constants.
This implies that $E$ is unbounded at point $A_4$.

In order to prove the unboundedness of the support, it will be more convenient to use  equation \eqref{U2}. Taking into account the properties of $E$ and $U_1$ established above, we obtain that near $A_4$ the right-hand side of this equation is equivalent to ${\rm const} \sqrt{|U_2-2w|}$. Since after separation of variables the integral over $U_2$ converges at the point $2w$, the support of the traveling wave is bounded.

Let us show that the derivatives of $U_1$ and $U_2$ are equal to zero at the end points of the support. For $U_2'$ this follows directly from the arguments of the previous section. From completely similar arguments it follows that the right-hand side of $U_1'$ is also equivalent to ${\rm const} \sqrt{|U_2-2w|}$ near $A_4$.

\begin{figure}[htb]
 \hspace{1.5cm}
\begin{minipage}{0.4\columnwidth}
\includegraphics[scale=0.6]{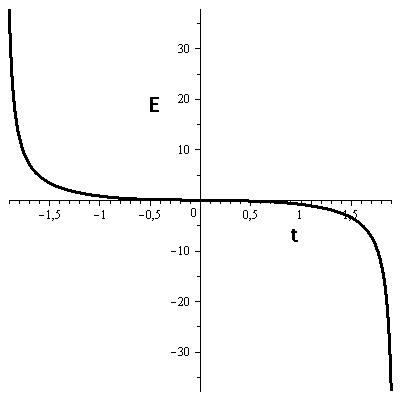}
\end{minipage}
\hspace{1.5cm}
\begin{minipage}{0.4\columnwidth}
\includegraphics[scale=0.6]{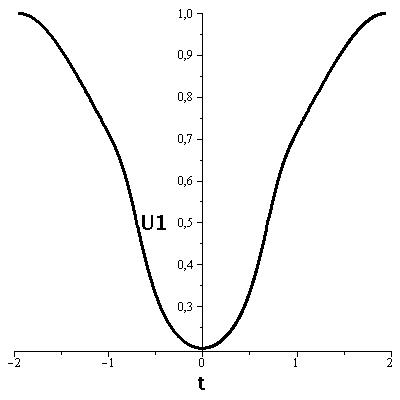}
\end{minipage}
\caption{Traveling wave for the second closure of parameter values from region 2, on the left  $E$, on the right  $U_1$. Everywhere $w=1$, $E(0)=0$, $U_1(0)=0.2$, $U_2(0)=0.7$. Traveling waves for the first closure do not exist in region 2}\label{2EU1}
\end{figure}

Fig. \ref{2EU1} illustrates the form of the traveling wave in region 2.

Case 3. Let us analyze the behavior of the trajectories near $A_3$. It is easy to see that \eqref{SpP} also holds near this point,
but now it follows that $U_1\sim w + C_1( U_2-w)$. Proceeding in a completely analogous manner to the previous points, we obtain from \eqref{EU2} that
near $A_3$ we have $E^2\sim \tilde C_1 U_2 +C_2$, $\tilde C_1(C_1), C_2$ are constants (from which follows the boundedness of $E$), and from \eqref{U2} the equivalence of the right-hand side of this equation ${{\rm const}}\cdot (U_2-w)^{-1}$, from which follows the boundedness of the carrier of the traveling wave. The unboundedness of the derivatives of all components of the solution at the end points of the support follows directly from \eqref{E}-\eqref{U2}.

\begin{figure}[htb]
 \hspace{1.5cm}
\begin{minipage}{0.4\columnwidth}
\includegraphics[scale=0.6]{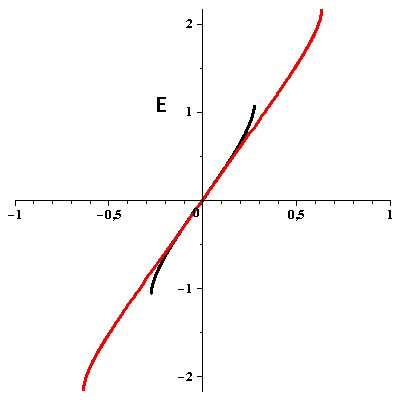}
\end{minipage}
\hspace{1.5cm}
\begin{minipage}{0.4\columnwidth}
\includegraphics[scale=0.6]{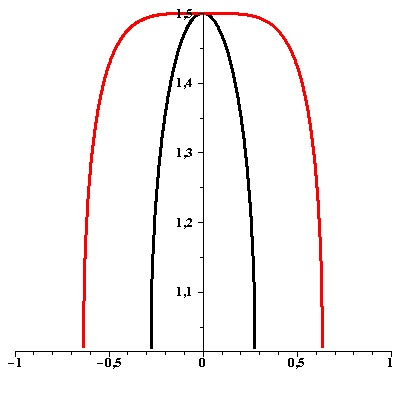}
\end{minipage}
\caption{Traveling wave for the second closure with parameters from domain 3 (red line), in comparison with the traveling wave for the first closure (black line), on the left $E$, on the right $U_1$. Everywhere $w=1$, $E(0)=0$, $U_1(0)=1.5$. For the first closure $U_2(0)=1.2$, for the second closure $U_2(0)=2$.}\label{3EU1}
\end{figure}

Fig.\ref{3EU1} illustrates  a traveling wave in region 3 in comparison with a traveling wave at $U_1=U_2$.

Thus, the proposition is proved. $\Box$

\medskip

Note that since the derivative $E'$ does not change sign for all the solutions constructed above, if we want to periodically continue the solution to the entire axis, then in all cases the component $E$ will be discontinuous, that is, according to \eqref{Kn}, the density component will contain a delta-shaped singularity. Since in region 3 the component $E$ increases, this singularity has a negative amplitude, and the resulting solution has no physical meaning. However, for domains 1 and 2 such solutions make sense, for them the components $U_1$ and $U_2$ turn out to be continuous on the entire axis, and for region 2 they are even smooth.

Thus, we see that for the closure at the second step the properties of the traveling waves change significantly in comparison with the original system of cold plasma dynamics, coinciding with the closure at the first step.

\section{Conclusion}\label{S7}
In this paper, we derive the equations of moment chains based on the kinetic equation of cold plasma and conduct its initial study for the case of one-dimensionality in space. For the closure at the first step, we obtain the standard hydrodynamic equations of cold plasma, the properties of whose solutions are well known. Our goal is to study the change in the properties of the solution upon closure at the next steps. We did the initial work to study the closure at the second step. Namely, we  studied two classes of solutions: linear in space and solutions of the traveling wave type, and we are convinced that the change in the properties of the solutions is very significant and cannot always be interpreted as an improvement in the smoothness properties of the solutions of the previous closure. On the other hand, the solutions we have considered are very specific (they  have infinite energy, and the traveling wave does not have initial smoothness) and cannot be considered as typical. Apparently, the question of whether the closure at the next step improves the properties of the solution can only be answered using numerical methods.

\section*{Acknowledgements}

 Supported by Russian Science Foundation  grant 23-11-00056 through RUDN University.


\begin{thebibliography}{99}
  \bibitem{AP56} { A.I. Akhiezer, R.V. Polovin},
{ Theory of wave motion of an electron plasma,}
JETP, 1956 \textbf {3} (5), 696.

\bibitem{ASR} D.Arsenio, L. Saint-Raymond, From the Vlasov-Maxwell-Boltzmann system to incompressible viscous electro-magneto-hydrodynamics -Vol.1, EMS Monographs in Mathematics, European Mathematical Society,  2019. DOI:10.4171/193

\bibitem{Bobylev}     A.V. Bobylev,  Boltzmann equation and hydrodynamics beyond Navier-Stokes,  Phil. Trans. R. Soc. A, 2018 {\bf 376} 20170227. DOI: 10.1515/9783110550986-008

    \bibitem{Cai}
Z. Cai, Y. Wang,  Regularized 13-moment equations for inverse power law models, Journal of Fluid Mechanics, 2020  {\bf 894} A12. DOI: 10.1017/jfm.2020.251

 \bibitem{CFL}
Z. Cai, Y. Fan, R. Li, Hyperbolic model reduction for kinetic equations. In: Recent Advances in Industrial and Applied Mathematics. SEMA SIMAI Springer Series, vol 1. Springer, Cham. 2022. DOI:10.1007/978-3-030-86236-7\, 8

\bibitem{Iord}
S.V. Iordanskii, On the Cauchy problem for the kinetic equation of plasma, Proceedings of the Steklov Mathematical Institute of the USSR, 1961 {\bf 60} 181-194.

 \bibitem{Koellermeier}
J. Koellermeier, R. P. Schaerer, M. Torrilhon, A framework for hyperbolic approximation of kinetic equations using quadrature-based projection methods, Kinetic and Related Models, 2014, {\bf 7}(3) 531-549. DOI:10.1137/14100110X

 \bibitem{Kuehn} C. Kuehn, {Moment closure - A brief review}, In: Control of Self-Organizing Nonlinear Systems, 253-271, Springer Cham, 2016. DOI:10.1007/978-3-319-28028-8\,13
\bibitem{Landau} L.D. Landau, On oscillations of electron plasma. JETP, 1946, {\bf 16} (7)
574-586.
\bibitem{Grad}
H. Grad. Asymptotic theory of the Boltzmann equation II Rarefied Gas Dynamics
(Proc. of the 3rd Intern. Sympos. Palais de l'UNESCO, Paris, 1962) Vol. I, 26-59.


\bibitem{Oraevskii}
V. Oraevskii, R. Chodura, W. Feneberg. Hydrodynamic equations for plasmas in strong
magnetic fields - I. Collisionless approximation, Plasma Physics, 1968 {\bf 10} 819-828.

\bibitem{RCh} O.S.Rozanova, E.V.Chizhonkov, On the conditions for the breaking of oscillations in a cold plasma,
ZAMP, 2021, {\bf 72}(1), 13. DOI:10.1007/s00033-020-01440-3
\bibitem{Rozanova:RChD_fric} O. S. Rozanova, E. V. Chizhonkov, M. I. Delova,  { Exact thresholds in the dynamics of cold plasma with electron-ion collisions},  AIP Conference Proceedings, 2021 {\bf 2302} (1), 060012. DOI:10.1063/5.0033619


\bibitem{SR}	L. Saint-Raymond, Hydrodynamic limits of the Boltzmann equation, Lecture Notes in Mathematics 1977,
Springer,	2009. DOI:0.1007/978-3-540-92847-8

\bibitem{Struchtrup}
H. Struchtrup, H. C. \"{O}ttinger, Thermodynamically admissible 13-moment equations. Physics of Fluids, 2022, {\bf 34} (1070-6631), 017105. DOI:10.1063/5.0078780

\bibitem{Torrilhon}
M. Torrilhon, H-Theorem for nonlinear regularized 13-moment equations in kinetic gas theory, Kinetic and Related Models, 2012, {\bf 5} (1) 185-201. DOI:10.3934/krm.2012.5.185

\bibitem{Kuz_Zakh} V. E. Zakharov, E. A. Kuznetsov, { Hamiltonian formalism for nonlinear waves}, 
    Phys. Usp., 1997, {\bf 40}(11)  1087-1116. DOI:10.1070/PU1997v040n11ABEH000304


\end{thebibliography}
\end{document}